\documentclass{lms}
\usepackage{amssymb,url,amsmath,amsfonts}
\usepackage[group-minimum-digits=4,group-separator=\text{,}]{siunitx}
\usepackage{bigstrut}
\usepackage{hyperref}
\usepackage{color}
\definecolor{mylinkcolor}{rgb}{0.5,0.0,0.0}
\definecolor{myurlcolor}{rgb}{0.0,0.0,0.75}
\hypersetup{colorlinks=true,urlcolor=myurlcolor,citecolor=myurlcolor,linkcolor=mylinkcolor,linktoc=page,breaklinks=true}

\newtheorem{theorem}{Theorem}
\newtheorem{computation}[theorem]{Computation}
\newtheorem{algorithm}[theorem]{Algorithm}

\def\CC{\mathbb{C}}
\def\FF{\mathbb{F}}
\def\PP{\mathbb{P}}
\def\ZZ{\mathbb{Z}}
\def\GL{\mathrm{GL}}
\def\PGL{\mathrm{PGL}}

\title[Zeta functions of quartic K3 surfaces over $\mathbb{F}_2$]
{A census of zeta functions of quartic K3 surfaces over $\mathbb{F}_2$}
\author{Kiran S. Kedlaya and Andrew V. Sutherland}
\date{April 29, 2016}
\classno{11M38, 14J28}

\extraline{Kedlaya was supported by NSF grant DMS-1501214, UCSD (Warschawski chair), and a Guggenheim Fellowship.  Sutherland was supported by NSF grants DMS-1115455 and DMS-1522526.}

\begin{document}

\maketitle

\begin{abstract}
We compute the complete set of candidates for the zeta function of a K3 surface over $\FF_2$ consistent with the Weil conjectures, as well as the complete set of zeta functions of smooth quartic surfaces over $\FF_2$. These sets differ substantially, but we do identify natural subsets which coincide. This gives some numerical evidence towards a
Honda-Tate theorem for transcendental zeta functions of K3 surfaces; such a result would refine a recent theorem of Taelman, in which one must allow an uncontrolled base field extension.
\end{abstract}

\section{Introduction and results}

For $X$ an algebraic variety over a finite field $\FF_q$ of characteristic $p$, the \emph{zeta function} of $X$ is the power series
\[
\zeta(X,T) = \exp \left( \sum_{n=1}^\infty \#X(\FF_{q^n}) \frac{T^n}{n} \right).
\]
A number of basic properties of $\zeta(X,T)$ are controlled by the now-proved \emph{Weil conjectures} (see for example \cite{osserman}); for example, $\zeta(X,T)$ always represents a rational function in $\mathbb{Q}(T)$.

Given a class of varieties over a particular field, it is natural to pose the \emph{inverse problem} asking which zeta functions consistent with the Weil conjectures actually occur. For abelian varieties of a given dimension $g$ over $\FF_q$, this question is resolved by  celebrated theorems of Tate \cite{tate} and Honda \cite{honda}: all such zeta functions occur provided that if $q \neq p$, one adds an extra condition on the factorizations over $\mathbb{Q}$ and $\mathbb{Q}_p$. (This condition always holds in the ordinary case; see \cite[\S 4]{haloui} for a concise statement of the condition and \cite{milne-waterhouse} for a thorough exposition.) By contrast, for curves of a given genus $g$, there are many additional constraints (the inequalities
$\#X(\FF_q) \geq 0$ and $\#X(\FF_{q^{mn}}) \geq \#X(\FF_{q^n})$, for example), and even the maximum value of $\#X(\FF_q)$ is unknown in most cases (see \cite{manypoints} for some results).

In this paper, we make a numerical investigation of the inverse problem for zeta functions of K3 surfaces over~$\FF_2$.
Recall that a \emph{K3 surface} over $\FF_q$ is a smooth, simply connected\footnote{Here we mean that the surface admits no nontrivial connected finite \'etale covers. This condition is needed only to eliminate abelian surfaces.} projective surface with trivial canonical bundle. The geometry of K3 surfaces is in many ways analogous to that of elliptic curves (they are Calabi-Yau varieties of dimensions 2 and 1, respectively). However, one key difference is that K3 surfaces cannot be uniformly described using a single geometric construction.
Instead, an infinite number of distinct constructions are required; we will focus mainly on the case of smooth quartic (degree 4) surfaces in $\PP^3$.

In the case of a K3 surface, one reads off from the Weil conjectures (plus properties of crystalline cohomology) the following constraints; see
\cite[Theorem~1]{taelman} for references and a sharper statement in case $q$ is not prime
(in the same vein as the Honda-Tate theorem).
\begin{theorem} \label{T:zeta K3}
Let $X$ be a K3 surface over $\FF_q$. Then $\zeta(X,T)$ has the form
\[
\frac{1}{(1-T)(1-qT)(1-q^2 T) q^{-1} L(qT)}
\]
for some polynomial $L(T) \in \ZZ[T]$ of degree $21$ with $L(0)= q$ having all roots on the unit circle.
\end{theorem}

One also has the following consequence of the Artin-Tate formula \cite{elsenhans-jahnel};
see Theorem~\ref{T:AT} for more discussion.

\begin{theorem} \label{T:ej}
With notation as in Theorem~\ref{T:zeta K3}, write $L(T) = (1-T)^r L_1(T)$ with $L_1(1) \neq 0$. Then $L_1(-1)$ is a perfect square (possibly $0$).
\end{theorem}

With these preliminaries in hand, we describe our computational results concerning zeta functions of K3 surfaces over $\FF_2$; the code used for these computations can be found in the repository \url{https://github.com/kedlaya/root-unitary}.
Our first computational result is an enumeration of Weil polynomials
based on a refinement of the search strategy described in~\cite{kedlaya-search}; see \S\ref{sec:tabulation} for details.
\begin{computation} \label{comp:abstract zeta}
The following sets are computed.
\begin{enumerate}
\item[(a)]
The set of polynomials $L(T)$ satisfying the conditions of 
Theorem~\ref{T:zeta K3} for $q=2$; it contains \num{2971182} elements.
\item[(b)]
The set of polynomials in (a) consistent with Theorem~\ref{T:ej};
it contains \num{2195801} elements.
\item[(c)]
The set of polynomials in (b) consistent with the inequalities  $\#X(\FF_q) \geq 0$ and $\#X(\FF_{q^{mn}}) \geq \#X(\FF_{q^n})$ (it suffices to impose the second condition for $(mn,n) \in \{(2,1), (3, 1), (4,2)\}$); it contains \num{1672565} elements.
\end{enumerate}
\end{computation}

Our second computational result is a lower bound for the inverse problem obtained by enumerating smooth quartic surfaces $X/\FF_2$ and computing $\zeta(X,T)$ directly by counting points in $X(\FF_{2^n})$;  see \S\ref{sec:point counting} for details.
\begin{computation} \label{comp:quartic zeta}
The following sets are computed.
\begin{enumerate}
\item[(a)]
The set of $\mathrm{PGL}_4(\FF_2)$-equivalence classes of smooth quartic surfaces over $\FF_2$;
it contains \num{528257} elements. 
\item[(b)]
The set of zeta functions of the surfaces in (a);
it contains \num{52755} elements and is a subset of the set found in Computation~\ref{comp:abstract zeta}(c).
\end{enumerate}
\end{computation}
With regard to (a), note that distinct $\mathrm{PGL}_4$-equivalence classes with the same zeta function may in fact give rise to isomorphic K3 surfaces: within the N\'eron-Severi lattice of a single K3 surface, the ample cone may contain multiple inequivalent divisors of degree 4. However, this cannot occur for $r=1$ (i.e., when $L'(1) \neq 0$). 
With regard to (b), note that in loose analogy with Tate's theorem that isogenous abelian varieties have the same zeta function \cite{tate}, a theorem of Lieblich and Olsson
\cite[Theorem~1.2]{lieblich-olsson} and Huybrechts \cite[Proposition~4.6]{huybrechts} states that K3 surfaces which are \emph{derived equivalent} (or \emph{Fourier-Mukai equivalent}) have the same zeta function.

Computation~\ref{comp:abstract zeta} provides a rich data set for investigating questions about zeta functions of K3 surfaces; for example, all possible values $1,\dots,10,\infty$ for the height of a K3 surface are realized by smooth quartics over $\FF_2$.
This said, the meaning of Computation~\ref{comp:quartic zeta} for the inverse problem for K3 surfaces over $\FF_2$ is unclear. However, it does yield some evidence towards a weaker form of the inverse problem suggested by Taelman \cite{taelman}. For $L(T)$ as in Theorem~\ref{T:zeta K3}, factor
$L(T) = \prod_i (1 - \alpha_i T)$ over $\CC$ and define the \emph{algebraic part} and \emph{transcendental part}
\[
L_{\mathrm{alg}}(T) 
:= \prod_{i: \alpha_i \in \mu_\infty} (1 - \alpha_i T),
\qquad
L_{\mathrm{trc}}(T)
:= \prod_{i: \alpha_i \notin \mu_\infty} (1 - \alpha_i T),
\]
where $\mu_\infty$ denotes the group of roots of unity in $\CC^\times$.
(By Computation~\ref{comp:abstract zeta}, for $q=2$ there are \num{73617} possible values of $L_{\mathrm{trc}}(T)$, whether or not we add conditions (b) and (c).) Then one can pose the inverse problem for $L_{\mathrm{trc}}(T)$ in place of $L(T)$, and in this case one has the following (conditional) partial solution
\cite[Theorem~2]{taelman}.
\begin{theorem} \label{T:taelman}
Assume that all K3 surfaces over finite extensions of $\mathbb{Q}_p$ have potential semistable reduction\footnote{This condition, which would be guaranteed if we were in equal characteristic 0 and is likely to hold in general, is only known for K3 surfaces with a polarization of small degree relative to $p$. See \cite[\S 2]{liedtke-matsumoto} for more discussion.}
in the sense of \cite[Definition~1]{taelman}.
Let $L_{\mathrm{trc}}(T) = \prod_i (1 - \alpha_i T)$ be a polynomial 
arising from some $L(T)$ as in Theorem~\ref{T:zeta K3}; if $q$ is not prime,
impose also the additional restrictions given in \cite[Theorem~1]{taelman}. Then for some positive integer~$n$,
the polynomial $\prod_i (1 - \alpha_i^n T)$ occurs as the transcendental part for some K3 surface over~$\FF_{q^n}$.
\end{theorem}

The proof of \cite{taelman} gives little insight as to whether the conclusion should necessarily hold with $n=1$. However, we can use our preceding computations to give a statement in this direction.
\begin{computation} \label{comp:comparison}
The following sets are computed.
\begin{enumerate}
\item[(a)]
The subset of Computation~\ref{comp:abstract zeta}(a) for which $L_{\mathrm{alg}}(T) = 1+T$, $L_{\mathrm{trc}}(1) = 2$, $L_{\mathrm{trc}}(-1) > 2$; it contains \num{1995} elements.
(Adding the conditions of Computation~\ref{comp:abstract zeta}(c) does not change this answer.)
\item[(b)]
The corresponding subset of Computation~\ref{comp:quartic zeta}(b); it contains the same \num{1995} elements.
\end{enumerate}
\end{computation}
The conditions imposed in Computation~\ref{comp:comparison} were chosen to partially (but not completely) eliminate the possibility that $L(T)$ arises from a K3 surface other than a smooth quartic by accounting for the Artin-Tate formula. 
It would be natural to continue the analysis by considering other families of K3 surfaces;
however, the Artin-Tate formula makes it difficult to produce enough examples to establish
that Theorem~\ref{T:taelman} always holds for $q=2$ with $n=1$. 
See \S\ref{sec:discussion} for elucidation of this point.

In another direction, one may hope to make similar calculations for $q>2$, but this poses significant technical challenges. Again, see \S\ref{sec:discussion} for further details.

\section{Tabulation of Weil polynomials}
\label{sec:tabulation}

Our tabulation of Weil polynomials broadly follows the search strategy described in \cite[\S 5]{kedlaya-search}; it is similar in spirit to the tabulation of number fields of prescribed signature, as in the work of Malle \cite{malle} and Voight \cite{voight}.
We briefly recall the strategy, indicate some new refinements which make it feasible to conduct much larger searches than previously possible, and discuss the computations performed.

\subsection{The search strategy}
This strategy attacks the problem of tabulating (not necessarily monic) integer polynomials
\[
P(T) = a_n T^n + \cdots + a_0
\]
of a fixed degree whose roots lie on the unit circle, with $a_0$ fixed and $a_1,\ldots,a_n$ constrained to lie in certain congruence classes (which could be a singleton set if the modulus is $0$, or all integers if the modulus is $1$).
All roots other than $\pm 1$ occur in conjugate pairs;
hence if $n$ is odd, then one of $\pm 1$ occurs with odd multiplicity, and we may reduce to the case where $n = 2m$ is even and $a_i = a_{n-i}$ for all $i$. 
(For example, for K3 surfaces we have $n=21$, so we first find reciprocal polynomials of degree 20, then multiply each one by $1+T$ and $1-T$ to generate the desired list.)
We may then write
\[
P(T) = T^m Q(T + 1/T), \qquad Q(T) = b_m T^m + \cdots + b_0
\]
for some integer polynomial $Q(T)$ with roots in the interval $[-2,2]$, with $b_m$ fixed and $b_0,\dots,b_{m-1}$ constrained to congruence classes $U_0,\ldots,U_{m-1}$.

To enumerate the set $S$ of polynomials $Q(T)$ we compute a tree with levels $0,\ldots,m$ in which each node $(b_m,\ldots,b_0)$ at level $m$ represents a polynomial $Q(T)=b_mT^m+\cdots+b_0$ in $S$, and nodes at level $i<m$ are labeled by tuples $(b_m,\ldots,b_{m-i})$ that are prefixes of their children $(b_m,\ldots, b_{m-i-1})$.
By Rolle's theorem, a necessary (but not sufficient) condition for the node $(b_m,\ldots,b_k)$ to have a descendant $(b_m,\ldots,b_0)$ at level $m$ is that the polynomial
\[
\sum_{j=0}^{m-k} \binom{k+j}{j} b_{k+j} T^j
\]
has all its roots in $[-2,2]$.
A general description of the algorithm appears below; a particular implementation is described in the next subsection.

\begin{algorithm} \label{search tree}
Given $b_m$ and congruence classes $U_0,\ldots,U_{m-1}$, enumerate the nodes of a tree with root $(b_m)$ at level $0$ in which each node at level $i>0$ is an integer tuple $(b_m,\dots,b_{m-i})$ with parent $(b_m,\ldots,b_{m-i+1})$ as follows.
\begin{enumerate}
\item[(a)]
Given a node $(b_m,\dots,b_{m-i})$, check whether the polynomial
\[
R(T) := \sum_{j=0}^{i} \binom{m-i+j}{j} b_{m-i+j} T^j
\]
 has all its roots in $[-2,2]$.
\item[(b)]
If the test in (a) passes and $i=m$, add $b_m T^m + \cdots + b_0$ to a list of return values.
\item[(c)]
If the test in (a) passes and $i < m$, compute an interval $I$ with the following property: 
for any values $b_{m-i-1},\dots,b_0$ such that $b_m T^m + \cdots + b_0$ has all roots in $[-2,2]$, one has $b_{m-i-1} \in I$.
Then take the children of this node to be the tuples $(b_m,\dots,b_{m-i-1})$
with $b_{m-i-1} \in U_{m-i-1} \cap I$ (this intersection may be the empty set).
\end{enumerate}
\end{algorithm}

\subsection{Implementation}

An implementation\footnote{While preparing this paper, we discovered a minor bug in the implementation accompanying \cite{kedlaya-search}. However, we corrected this bug in the current implementation and reconfirmed all of the computational results.}
 of the aforementioned search strategy (available for download) is described in \cite{kedlaya-search}. The implementation we use here differs from the prior one in several theoretical and practical aspects, which we now discuss.
\medskip

\begin{itemize}
\item
In \cite{kedlaya-search}, the interval $I$ in Algorithm~\ref{search tree}(c) is constructed using linear and quadratic inequalities on the power sums $s_1,\dots,s_{i+1}$, as computed from $b_m,\dots,b_{m-i-1}$ via the Newton identities; 
note that $m$ appears explicitly in the identities, so these inequalities typically carry more information than simply requiring that the one-step extension conform to Rolle's theorem. See
\cite[\S 5]{kedlaya-search} for the precise list of inequalities used. 
\smallskip

We add a new condition which is linear in the $b_j$ rather than the $s_j$: if $Q(T)$ has all roots in $[-2,2]$, then $b_m Q(2+T)$ and $(-1)^m b_m Q(2-T)$ have all coefficients nonnegative (consistent with Descartes's rule of signs).
\medskip

\item
In \cite{kedlaya-search}, the test in Algorithm~\ref{search tree}(a) is conducted as follows. Let $S(T)$ be the polynomial obtained from $R(T)$ by removing all multiple roots and all factors of $T \pm 2$.
Form a Sturm sequence $S_0,S_1,\dots$ as follows:  $S_0 := S$, $S_1 := S'$, and for $j>1$, let $S_i$ be a negative scalar multiple of the remainder of $S_{j-2}$ modulo $S_{j-1}$, stopping just before appending the zero polynomial to the sequence.
We then invoke Sturm's theorem  \cite[Theorem~2.50]{basu-pollack-roy}: the number of roots of $S(T)$ in $[-2,2]$ is the difference between the numbers of sign changes in the sequences $S_0(-2), S_1(-2), \dots$ and $S_0(2) ,S_1(2),\ldots$.
\smallskip

In the current implementation, we note that when testing $R$, we have the prior information that $R'(T)$ has all roots in $[-2,2]$ and that $R(2)$, $R(-2)$ have the correct signs (thanks to the previous point). We thus need only test that $R$ has all real roots; this eliminates some polynomial evaluations at $\pm 2$. We may also take $S = R$, since the fact that Sturm's theorem remains valid as a count of roots without multiplicity; this avoids duplication of the Euclidean algorithm. 
We next observe that $R$ has all real roots if and only if for $j=1,\dots,i$,
$\deg(S_j) = i-j$ and the leading coefficients of $S_0$ and $S_j$ have the same sign; this allows for an early abort. Finally, we note that if the early abort happens due to a sign discrepancy (rather than a degree discrepancy) at $S_j$,
then this discrepancy does not depend on $b_k$ for 
$k \leq m-i-1 + \max\{0, i-2j+2\}$; if $i-2j+2 > 0$, we may thus back up the tree traversal and discard all nodes below $(b_m,\dots,b_{m-2j+2})$ without losing any of the return values. (If $b_m > 0$ and nodes are traversed in lexicographic order, then one may replace $i-2j+2$ with $i-2j+1$ and $b_{m-2j+2}$ with $b_{m-2j+1}$ in the previous analysis.)
\medskip
\item
In \cite{kedlaya-search}, the implementation consisted of an interpreted component in \texttt{Sage} \cite{sage} performing high-level user interaction and
a compiled component in \texttt{Cython} \cite{cython} for mid-level computations.
Some low-level computations, such as Sturm sequences, were farmed out to compiled components of the \texttt{Sage} library, notably \texttt{PARI} \cite{pari}.\footnote{The PARI/GP project includes both the C library \texttt{PARI} and the interpreted \texttt{GP} language. \texttt{Sage} interfaces directly with the C library.}
\smallskip

In the current implementation, we incorporate a third component, written in C using the \texttt{FLINT} library \cite{flint}. This component absorbs most of the work of the \texttt{Cython} layer (which remains to provide wrappers around the C code) and completely supplants the use of \texttt{PARI}.

\medskip
\item
In \cite{kedlaya-search}, parallelization via work-stealing is suggested but not implemented; we provide this in the current implementation.
Given a pool of threads, we initially assign the entire search tree to one thread, then iterate the following steps until no active threads remain.
\begin{enumerate}
\item[(i)]
In parallel,
each thread which is active (i.e., has been assigned a subtree of the search tree)
performs a depth-first search to find one polynomial in its remaining search space, going inactive if none exist.
\item[(ii)]
In serial, each inactive thread solicits work from a randomly chosen active thread by removing a branch (as close to the root as possible) from the latter's subtree.
\end{enumerate}
\end{itemize}

\subsection{Computations}
\label{subsec:computations}

We now describe in detail some computational results obtained using this search strategy.
The reported computations were carried out on a 24 core Intel Xeon X5690 3.47GHz machine with 192GB of memory. The parallel implementation was run using 512 threads; this provided an 8--10$\times$ speedup.

Computation~\ref{comp:abstract zeta} was completed in
under 1 hour. In addition to the \num{1485591} polynomials of degree 20 that were found, the search tree found an additional \num{2149281061} leaves at smaller depths; that is, the dead ends outnumber the solutions by a factor of nearly \num{1500}. This suggests that there may still be substantial room to optimize the choice of the intervals in Algorithm~\ref{search tree}(b) (see below).

As an additional test of the implementation (which helped expose some bugs during development), we computed the set of monic polynomials $L(T) \in \ZZ[T]$ with degree in $\{1,3,\dots,21\}$ with all roots on the unit circle; it contains $\num{78670}$ elements.
By Kronecker's theorem, these polynomials factor as products of cyclotomic polynomials; it is thus easy to confirm this answer using an independent script that forms these products directly.

\subsection{Possible refinements}

The discussion above suggests that there remains substantial room for improvement in the 
choice of intervals in Algorithm~\ref{search tree}(b). As noted in \cite{kedlaya-search}, similar searches in other contexts often make use of linear programming methods; we have not investigated this direction.

There is also a question of whether Sturm sequences are the optimal method for the test in Algorithm~\ref{search tree}(a). For one, Sturm sequences involve multiprecision integers, in part due to the systemic appearance of certain large powers. In principle, these powers can be removed explicitly, thus reducing the computational complexity \cite[Algorithm~3.3.1]{cohen}; in practice, we find (in this particular setting) that the Gaussian content is substantially larger than predicted by general arguments, so we prefer to compute it explicitly. (One could try mixing the two approaches, but in \cite[\S 3.3]{cohen} it is suggested that this gives inferior results.)

More seriously, there is a question as to whether Sturm sequences are superior to root isolation methods based on the Budan-Fourier theorem \cite[Theorem~2.35]{basu-pollack-roy} (see also \cite[\S 10.4]{basu-pollack-roy}). We have chosen Sturm sequences in part for ease of implementation, but also because they are better suited to the task at hand. To wit,
root isolation methods can easily generate certificates that guarantee the existence of certain real roots (using sign changes), but have more difficulty generating certificates that guarantee the failure of a polynomial to have all of its roots in an interval. By contrast, Sturm's theorem provides certificates of the latter type easily using the early-abort mechanism described above. That said, it may be that a well-crafted strategy using root isolation (e.g., one which uses the positions of the roots of $R'$ to help isolate the roots of~$R$) would work better in the long run.

In order to get a fuller parallel speedup, some refinement of the parallelization mechanism  is needed. For example, we currently only interrupt a process when it finds a solution or exhausts its search space; this is suboptimal in certain use cases where the search tree is very large but the number of solutions to be found is small.

\section{Point counting}
\label{sec:point counting}

We computed the zeta function of every K3 surface over $\FF_2$ that arises as a smooth quartic surface $X$ in $\PP^3$ by counting points on these surfaces over extension fields $\FF_{2^n}$, with~$n$ ranging over a set of values sufficient to uniquely determine the zeta function $\zeta(X,T)$ given by Theorem~\ref{T:zeta K3} and Theorem~\ref{T:ej} (up to $n=19$ in the worst case).
This computational problem naturally breaks down into tasks: (1) enumerate a complete set $S$ of smooth surfaces defined by homogeneous quartic polynomials $f\in \FF_2[w,x,y,z]$, up to $\PGL_4$-equivalence; (2) compute $\#X(\FF_{2^n})$ for $X\in S$ and suitable values of~$n$.

\subsection{Determining $\PGL_4$-orbits of homogeneous quartics}

There are $\binom{7}{3}=35$ homogeneous quartic monomials $w^ax^by^cz^d$, one for each quadruple of nonnegative integers $(a,b,c,d)$ with $a+b+c+d=4$.
If we order the quadruples $(a,b,c,d)$ lexicographically, each homogeneous quartic $f\in\FF_2[w,x,y,z]$ can be uniquely identified with a bit-vector $v:=v(f)\in \FF_2^{35}$ indexed by quadruples $(a,b,c,d)$ for which $f=\sum v_{(a,b,c,d)}w^ax^by^cz^d$;
the bit-vector $v$ can be conveniently encoded as an integer in $[0,2^{35})$ and we order them accordingly (this is just the lexicographic order on $\{0,1\}^{35}$).

The group $\PGL_4(\FF_2)$ acts on the set of homogeneous quartics $f(w,x,y,z)$ via linear change of variables.
In terms of the corresponding set $V:=\FF_2^{35}$ of vectors $v(f)$, each element of $\PGL_4(\FF_2)$ corresponds to an invertible linear transformation of $V$ that can be explicitly represented as an invertible $35\times 35$ matrix.
We may thus identify $\PGL_4(\FF_2)$ with a subgroup~$G$ of $\GL_{35}(\FF_2)$ of order $2^6(2^2-1)(2^3-1)(2^4-1)= \num{20160}$.
As we are only interested in the quartic surfaces $f(w,x,y,z)=0$ up to isomorphism, it suffices to consider the $G$-orbits of $V$, each of which may be uniquely represented by a lexicographically minimal $v$.
The number of $G$-orbits can be computed via Burnside's lemma as
\begin{equation}\label{eq:burnside}
\#(V/G) = \frac{1}{\#G}\sum_g \#V^g = \frac{\#C}{\#G}\sum_C (\#\FF_2)^{\dim_1(C)} = \num{1732564},
\end{equation}
where the first sum is over group elements, the second sum is over conjugacy classes, and $\dim_1(C)$ denotes the dimension of the $1$-eigenspace of the conjugacy class $C$.  There are only 14 conjugacy classes in $\PGL_4(\FF_2)$, so the second sum is trivial to compute.

To find lexicographically minimal representatives for each orbit we simply enumerated every orbit using a bitmap with $2^{35}$ entries; this took less than 2 days.
We note that this brute-force approach is not feasible for finite fields larger than $\FF_2$.
Indeed, determining a set of unique orbit representatives is already a nontrivial problem over $\FF_3$ (the vector space $\FF_3^{35}$ contains $3^{35}\approx 2^{55.5}$ elements).  However, determining the cardinality of this set via \eqref{eq:burnside} is quite feasible for values of $q>2$; for example, over $\FF_3$ there are \num{4127971480} orbits, and over $\FF_5$ there are \num{100304466278983}. 

Having compiled a complete list of $\PGL_4(\FF_2)$-orbits of homogeneous quartics, we then want to restrict to those that define a K3 surface; this amounts to discarding orbits represented by a vector $v(f)$ for which the polynomial $f\in \FF_2[w,x,y,z]$ is not irreducible, or for which the singular locus defined by the Jacobian matrix of $f$ is nonempty (the latter implies the former but the former is often easier to check).
These conditions are straightforward to apply, and we quickly find that \num{528257} of the \num{1732564} orbits satisfy them; these constitute our set~$S$ of smooth plane quartic surfaces $X/\FF_2$ in $\PP^3$.

\subsection{Counting points on quartic surfaces over $\FF_2$}

Given a smooth quartic surface $X\in S$ defined by $f(w,x,y,z)=0$, our basic strategy for computing $\#X(\FF_{2^n})$ is elementary: iterate over pairs $(x_0,y_0)\in \FF_{2^n}^2$ and for each pair determine the number of roots of the polynomial $g(w):=f(w,x_0,y_0,1)=g(y)\in \FF_{2^n}[w]$ that lie in $\FF_{2^n}$ (of course we also need to account for points with $z=0$, but this reduces to the much easier problem of counting points on a curve in $\PP^2$ and takes negligible time).
To count the roots of $g(w)$ we use Zinoviev's formulas \cite{zinoviev}, which for low degree polynomials $g$ over $\FF_{2^n}$ gives explicit $n\times n$ systems of linear equations over $\FF_2$ whose solutions correspond to the roots of $g$, based on Berlekamp's algorithm for factoring polynomials over finite fields of small characteristic using linear algebra~\cite{berlekamp}.

One might generically expect $g$ to have degree 4, but in fact this is not the case.
For all but 34 of the surfaces in $S$, the degree of the defining polynomial $f(w,x,y,z)$ in $w$ is at most 3 (note that our lexicographic ordering minimizes the degree in $w$).
In the typical case where $g(w)$ is a cubic, after making it monic and applying a linear change of variable we may assume $g(w)=w^3+g_1w+g_0$.
It is then feasible to precompute a lookup table $T$ indexed by pairs $(g_0,g_1)\in \FF_{2^n}^2$ whose entries record the number of roots of $w^3+g_1w+g_0$ in $\FF_{2^n}$.
Each entry in $T$ is an integer in $[0,3]$ that can be encoded in 2 bits, thus the total size of $T$ is $2^{2r+1}$ bits; even for $r=19$, this is reasonably small (64 GB).
The time to compute $T$ is actually less than the time to instantiate $f(w,x_0,y_0,1)$ at every pair $(x_0,y_0)\in \FF_{2^n}^2$; we can accelerate this computation by enumerating the pairs $(g_0,g_1)$ in an order that makes it convenient to compute the matrices appearing in Zinoviev's formulas.
Thus makes the computation of $T$ worthwhile even for a single surface $X$, and we can reuse the same table $T$ for every $X\in S$.
The cost of point counting is then dominated by the time to evaluate $f(w,x_0,y_0,1)$.

Given that we wish to count points on a fairly large set of surfaces (all of $S$ for $n\le 12$ and subsets $S_n\subseteq S$ for $n>12$), rather than iterating over surfaces $f(w,x,y,z)=0$ and counting points on each, which involves iterating pairs $(x_0,y_0)\in \FF_{2^n}^2$ and evaluating $f(w,x_0,y_0,1)$, we reverse the order of iteration and loop over pairs $(x_0,y_0)$ and for each pair count the solutions to $f(w,x_0,y_0,1)=0$ over $\FF_{2^n}$ for every surface $f(w,x,y,z)=0$ in our set, keeping a running total of points for each surface as we go.
This allows us to instantiate the 35 homogeneous quartic monomials at $x=x_0,y=y_0,z=1$ just once for each pair $(x_0,y_0)$, and then for each polynomial $f(w,x,y,z)$ compute $f(w,x_0,y_0,1)$ as an $\FF_2$-linear combination of these, equivalently, as a sum of a subset of them, which is very fast.

This algorithm is trivially parallelizable (with linear speedup), and running on 32 cores it takes only two days to compute $\#X(\FF_{2^n})$ for $1\le n\le 12$ and all \num{528257} surfaces $X\in S$. From these point counts, for each $X \in S$ we can write
$L(T) = 1 + a_1 T + \cdots + a_{21} T^{21}$ with $a_1,\dots,a_{12}$ known.
In most cases, the existence of a sign $\epsilon \in \{+, -\}$ such that
$a_{21-i} = \epsilon a_i$ then determines $L(T)$ uniquely;
the exceptions are the cases where $a_{9} = a_{10} = 0$. 
For $n=13,\dots,21$, let $S'_n$ be the subset of $S$ consisting of exceptions
for which $a_{21-n} \neq 0$ and $a_{20-n} = \dots = a_{10} = 0$;
let $S_n$ be the subset of $S_n$ for which both choices for $\epsilon$ give polynomials compliant with Theorem~\ref{T:zeta K3}, Theorem~\ref{T:ej}, and Theorem~\ref{T:AT} (which for degree 4 means that if $r=1$, then $L_1(1)$ is a square).
The sizes of these sets are listed in Table~\ref{table:misses}.

For $n=13,\dots,19$, we reran the previous algorithm to compute $\#X(\FF_{2^n})$ for each $X \in S'_n$; the most time-consuming computation was for $S'_{19}$, which took about six days running on a machine with 32 Intel Xeon E5-2687Wv2 3.4GHz cores and 256GB of memory.
(It would have been sufficient to consider $x \in S_n$, but the extra computations serve as a consistency check.)
Note that the $\num{1876}$ tuples $(\#X(\FF_2), \dots, \#X(\FF_{2^{12}}))$ represented by $S_{13} \cup \cdots \cup S_{19}$ only give rise to $\num{2071}$ different zeta functions; that is, in the vast majority of these cases only one of the two sign choices is realized. This suggests that there may be further theoretical restrictions on zeta functions that we have not yet taken into account (e.g., interaction between the Newton polygon and the order of the Brauer group).

\begin{table}[ht]
\begin{tabular}{c|ccccccccc}
$n$ & 13 & 14 & 15 & 16 & 17 & 18 & 19 & 20 & 21 \\
\hline\bigstrut[t]
$\#S'_n$ & \num{38225} & \num{16555} & \num{8281} & \num{3608} & \num{2011} & \num{857} & \num{283} & \num{0} & \num{96} \\
$\#S_n$ & \num{17795} & \num{7315} & \num{3611} & \num{1435} & \num{1016} & \num{470} & \num{125} & \num{0} & \num{0}\\
\hline
\end{tabular}
\caption{Cardinalities of sets $S'_n, S_n$}
 \label{table:misses}
\end{table}

We should also mention an important low-level optimization to speed up arithmetic in $\FF_{2^n}$ that we used: the ``carry-less multiplication" instruction PCLMULQDQ now available on Intel processors (since 2010) speeds up multiplication in $\FF_{2^n}$ quite dramatically (by a factor of up to 10 for the values of $n$ that we used).

\section{Further discussion}
\label{sec:discussion}

In light of the preceding computations, we resume the discussion from the introduction
concerning the inverse problem for zeta functions of K3 surfaces.

As noted earlier, smooth quartics give rise to only one out of infinitely many algebraic families of K3 surfaces. This is due to the fact that, just as for abelian varieties, in order
to represent the moduli problem one must consider \emph{polarized} K3 surfaces.
The \emph{degree} of the polarization equals its self-pairing in the N\'eron-Severi lattice; each value of the degree corresponds to a single irreducible component of the moduli space of K3 surfaces. For the first few degrees, the generic polarized K3 surface of that degree can be described as follows.
(There are also some nongeneric possibilities; see below for some discussion of the case of degree 4.)
\begin{itemize}
\item
Degree 0: an elliptic K3 surface.
\item
Degree 2: a double cover of $\PP^2$ (or a twist) branched over a smooth sextic (degree 6) curve.
\item
Degree 4: a smooth quartic in $\PP^3$.
\item
Degree 6: a smooth transversal intersection of a quadric and a cubic in $\PP^4$.
\item
Degree 8: a smooth transversal intersection of three quadrics in $\PP^5$.
\end{itemize} 

For $X$ a K3 surface over $\FF_q$, define $L(T)$ as in Theorem~\ref{T:zeta K3}
and $r, L_1(T)$ as in Theorem~\ref{T:ej}. The interaction between these invariants and the degree is governed by the Artin-Tate formula, which was used already to deduce
Theorem~\ref{T:ej} (by comparing \eqref{eq:AT} over $\FF_q$ and $\FF_{q^2}$);
here is the full statement, as in \cite[Proposition~6]{elsenhans-jahnel}.

\begin{theorem}[(Artin-Tate formula)] \label{T:AT}
If $X$ satisfies the Tate conjecture,\footnote{This is known to hold except in certain cases in characteristic 2; see \cite{madapusi}. The case of characteristic 2 is apparently resolved by a very recent preprint \cite{kim-madapusi}. In any case,
the formulation in \cite[Proposition~6]{elsenhans-jahnel} is made carefully so as not to require the Tate conjecture.} then 
\begin{equation} \label{eq:AT}
L_1(1) = |\Delta| \#\mathrm{Br}(X),
\end{equation}
where $\Delta$ denotes the discriminant of the N\'eron-Severi lattice and $\mathrm{Br}(X)$ denotes the Brauer group; the latter is finite and its order is a perfect square. 
Also, the rank of the N\'eron-Severi lattice equals $r$; in particular, if $r=1$
then $X$ admits a unique polarization, and its degree is equal to $\Delta$.
\end{theorem}

We can now justify the choice of the conditions in Computation~\ref{comp:comparison}:
besides smooth quartics, some other sources of degree-4 K3 surfaces include desingularization of singular quartics with only isolated rational singularities, which all have $r>1$; and 
double covers of quadrics branched along $(4,4)$ curves, which all have either $r>1$, $L_{\mathrm{alg}}(T) \neq 1+T$, or
$\Delta = -4$.

We now justify our previous assertion that Theorem~\ref{T:AT} makes it difficult to produce enough examples to establish that Theorem~\ref{T:taelman} always holds for $q=2$ with $n=1$. 
Among the possibilities for $L_{\mathrm{trc}}(T)$ allowed by Computation~\ref{comp:abstract zeta}(c), there exist cases where
\[
\deg L_{\mathrm{trc}}(T) = 20, \qquad L_{\mathrm{trc}}(1) \in\{307, 367, 463\},   \qquad L_{\mathrm{trc}}(-1) = 3.
\]
By Theorem~\ref{T:ej}, these can occur only with $L_{\mathrm{alg}}(T) = 1+T$. Since $307$, $367$, and $463$ are prime, these would have to occur for K3 surfaces of degrees $2 \times 307$, $2 \times 367$, and $2 \times 463$, respectively; however, the moduli spaces of polarized K3 surfaces of these degrees are of general type  \cite[Theorem~1]{gritsenko}, so constructing explicit points on them may be difficult. A more promising approach would be to make explicit the constructions used in \cite{taelman}, which involve lifting to characteristic 0, as in the proof of Honda's theorem.


Finally, we discuss the prospects for making similar calculations for $q>2$.
On the side of enumerating candidate zeta functions, there seems to be a bit of room to enlarge $q$; for example, the following computation for $q=3$ took 2.5 days (wall time, parallelized as in \S\ref{subsec:computations}).

\begin{computation}
The following sets are computed.
\begin{enumerate}
\item[(a)]
The set of polynomials $L(T)$ satisfying the conditions of 
Theorem~\ref{T:zeta K3} for $q=3$;
it contains $\num{75936610}$ elements, representing $\num{6867811}$ distinct values of $L_{\mathrm{trc}}(T)$.
\item[(b)]
The set of polynomials in (a) consistent with Theorem~\ref{T:ej};
it contains $\num{52980075}$ elements.
\item[(c)]
The set of polynomials in (b) consistent with the inequalities  $\#X(\FF_q) \geq 0$ and $\#X(\FF_{q^{mn}}) \geq \#X(\FF_{q^n})$ (it suffices to impose the second condition for $(mn,n) = (2,1)$); it contains $\num{49645728}$ elements.
\end{enumerate}
\end{computation}

On the side of computing zeta functions of quartics, we must emphasize that the (optimized) brute-force approach we used to compute the zeta functions of all quartic surface over $\FF_q$ by point counting is not feasible for $q>2$, particularly in the exceptional cases where one must go beyond $\FF_{q^{11}}$ to resolve the sign ambiguity.
For $q=3$ there are already more than 2000 times as many $\PGL_4$-orbits of quartics to consider, and the time to compute $\#X(\FF_{q^{11}})$ will be larger by a factor of at least 100 (more than 2000 for $\FF_{q^{19}}$).
One should instead look to methods based on $p$-adic cohomology, as in~\cite{akr}; these have recently been made practical\footnote{The implementation used in \cite{costa} does not allow $q=2$, so we cannot use it to check Computation~\ref{comp:quartic zeta}.} for K3 surfaces \cite{costa}. 

\acknowledgments{This work was carried out at ICERM during the fall 2015 semester program ``Computational aspects of the Langlands program''. We thank Edgar Costa, David Harvey, Brendan Hassett, Christian Liedtke, Sebastian Pancratz, Matthias Sch\"utt, Lenny Taelman, and Yuri Tschinkel for helpful discussions.}

\affiliationone{
   Kiran S. Kedlaya\\
   Department of Mathematics \\
   University of California, San Diego \\
   9500 Gilman Drive \#0112 \\
   La Jolla, CA\ 92093\\
   United States of America \\
   \email{kedlaya@ucsd.edu}}
\affiliationtwo{
    Andrew V. Sutherland\\
    Department of Mathematics\\
    Massachusetts Institute of Technology\\
    77 Massachusetts Avenue\\
    Cambridge, MA \ 02139\\
    United States of America\\
    \email{drew@math.mit.edu}}

\end{document}